# How Disciplinary Norms Influence Mathematicians' Views of Programming in Undergraduate Mathematics

Jan-Fredrik Olsen[1,2] (jan-fredrik.olsen@math.lu.se)

Tor Ole B Odden[1] (t.o.b.odden@fys.uio.no)

*Programming is deeply embedded in contemporary mathematical practice, yet its epistemic status in university mathematics teaching remains contested. Little is known about how mathematicians themselves understand the legitimacy of programming in their professional work, and how these views shape their teaching. We address this gap through semi-structured interviews with 15 mathematicians at a Northern European university with over two decades of systematic integration of programming across STEM subjects. Drawing on Cultural-Historical Activity Theory and Communities of Practice, we examine how mathematicians articulate the role of programming across research and teaching.*

*We identify four epistemic archetypes – classical pure, classical applied, computational applied, and computational pure – each expressing coherent norms governing legitimate use of programming. Across archetypes, light-touch programming (e.g., numerical exploration) was widely used in research but largely invisible in teaching, where legitimacy was tied to more "substantive" integration. We argue that this gap reflects epistemic continuity across practices combined with teaching's focus on established analytic outcomes, which reduces the epistemic visibility of computational work.*

*Given how mathematicians articulate legitimate uses of programming, our findings suggest that integration is most widely accepted when substantive programming is taught in dedicated courses, while light epistemic use – such as numerical exploration – is used to support learning in traditional theoretical courses. More extensive computational work is generally viewed as fitting naturally within specialised numerical or computational mathematics courses.*

## Introduction

The notion of programming has shaped mathematical practice for nearly two centuries. Already in 1843, Ada Lovelace argued that Babbage's Analytical Engine could "*weave algebraical patterns just as the Jacquard loom weaves flowers and leaves*" by following sequences of symbolic instructions (as reprinted in Lapham's Quarterly; Lovelace, n.d.). She anticipated that such programmable mechanisms would render mathematical truths "*more speedy and accurate [in] practical application.*"

The legitimacy of computer-aided computation in mathematics draws both on their value for real-world and theoretically oriented mathematicians. Indeed, such computations are essential to mediate between mathematical

---

[1] Center for Computing in Science Education, University of Oslo, Sem Saelands vei 24, Fysikkbygningen, 0371 Oslo, Norway
[2] Centre for Mathematical Sciences, Lund University, Box 118, 22100 Lund, Sweden

models and real-world data, and afford opportunities to explore theoretical conjectures, visualise structures, and test ideas that may seem analytically intractable (Broley et al, 2018; Lockwood et al., 2019).

Yet the role of computations in mathematics is controversial. Paul Halmos famously reacted in disgust to Appel and Haken's computational proof of the Four Color Theorem by exclaiming "*The present proof relies in effect on an Oracle, and I say down with Oracles! They are not mathematics*" (Hersh, 1997, p. 54). As Hersh explains, Halmos' judgement is not based on doubts about correctness but from concerns about what counts as a '*proper*' mathematical argument. This is in line with findings that mathematicians evaluate use of computations not only on logical grounds, but also on aesthetic and epistemic ones – including transparency, explanatory depth, beauty, and purity (Burton, 2004; Hamami & Morris, 2020).

These tensions appear to carry over to postsecondary education, and seem to be even more pronounced there; existing research shows that engagement with programming is even lower in teaching than research, especially within theoretically oriented courses (Buteau et al., 2014; Broley et al., 2018). Understanding this discrepancy is important if programming is to be integrated into mathematics education in a meaningful and sustainable way. However, the connection between the research and teaching practices of mathematicians remains poorly understood, with prior work suggesting that mathematicians may enact different epistemologies across these domains (see, e.g., Burton, 2004).

In this paper, we interview 15 mathematicians to investigate how mathematicians articulate the legitimacy of programming across their research and teaching practices. Our study is situated at the University of Oslo (UiO), where programming and scientific computation have been formally integrated into mathematics programmes for more than twenty years (Malthe-Sørenssen et al., 2016). In their call for research on the use of programming in undergraduate mathematics education, Lockwood et al. (2021) explicitly suggest this as a critical site that "offers a fascinating look at program-level cooperation that can be studied more systematically."

In our analysis, we draw on a socio-cultural theoretical framework to construct an epistemic matrix that characterises mathematicians by their dominant object orientation (what they study) and tool orientation (how they study it). This matrix helps us reduce heterogeneity in the data and to identify epistemic rules that shape how programming is legitimised across research and teaching. We use this framework to trace a set of recurring tensions – such as the persistence of research-based norms in teaching and the diminished epistemic visibility of programming when mathematical activity shifts from producing knowledge to teaching established outcomes. The research questions guiding this study are presented at the end of the Theoretical framework section.

# Background

## Mathematicians' practices in mathematics education research

Mathematicians' practices form an established line of inquiry in mathematics education research, with important early contributions by Burton (2004) and Nardi (2008). Building on Moschkovich's (2007) notion of professional discourse practices, Rasmussen et al. (2015) defined disciplinary practices as "the ways in which mathematicians go about their profession" (p. 264). For an overview of research in this area, see Weber et al. (2020).

Burton (2004) found that mathematicians' research practices were far more heterogeneous and exploratory than their classroom practices suggest. In her analysis of seventy interviews, she reports that mathematicians describe their research in terms of exploration, intuition, collaboration, conjecturing, and invention, yet adopt a transmissive model when teaching – a model that presents mathematics as logical structures in a Platonic landscape "missing entirely the stumbling human process that created those results in the first place" (Smith & Hungwe, 1998, as cited in Burton, 2004, p. 87). Anticipating our theoretical framework, we refer to this as the *object-outcome shift* of mathematics education.

Research shows that mathematicians' self-described epistemologies (see, e.g., Polya, 2014) capture only part of a far more diverse set of practices, beliefs, and norms. Burton (2004) conceptualises mathematical knowing through a broad epistemological model that incorporates social relations, aesthetics, and recognition of differing approaches to inquiry. In contrast, Weber et al. (2014) offer an epistemic model centred on epistemic aims and actions – what mathematicians actually do to gain conviction. Their findings show that mathematicians sometimes rely on empirical or authoritative evidence and sometimes fail to achieve full conviction even from proofs, underscoring the complexity of their epistemic behaviour.

Building on this work, Weber et al. (2020) highlight several methodological challenges for studying mathematicians' practices. Beyond the heterogeneity problem already noted, they identify the accuracy problem (mathematicians may unwittingly misrepresent their own practices), the identification problem (what counts as a 'mathematician' varies across contexts), the advanced-content problem (specialised language can obscure meaning for outsiders), and the interdisciplinarity problem (adequately analysing practice requires engagement with history and philosophy of mathematics). These challenges collectively suggest that research in this area benefits from close collaboration with professional mathematicians (see also Darragh, 2022).

## Programming in mathematicians' practices

Despite increasing attention to mathematicians' practices in mathematics education research, only a small body of work has examined their use of programming. This gap is noteworthy given that the epistemic status of computation in mathematics has long been contested, as the Four Color Theorem controversy alluded to in the introduction shows.

Empirical evidence suggests that this tension extends to mathematicians' everyday use of programming. In a survey of 302 Canadian mathematicians, Buteau et al. (2018) found that 18% used programming in their teaching and 43% in their research. In follow-up interviews, Broley et al. (2018) reported that programming was viewed as having little place in pure mathematics courses, while being more accepted in service courses or applied courses. The authors attribute this to a 'status gap', shaped by institutional and cultural factors. Together, these studies point to systematic differences both in the kinds of mathematical objects mathematicians prioritize and the tools they view as legitimate.

At the same time, studies show that mathematicians who do use programming emphasise both pragmatic and epistemic benefits. Broley et al. (2018) and Lockwood et al. (2019) report that mathematicians use programming to experiment, check conjectures, produce visualisation and make approximative computations. Broley et al. (2018) also highlight perceived pedagogical advantages: by opening up black boxes and building their own tools, students can "*feel a sense of empowerment and excitement that may further enhance their experience*".

## Programming in mathematics education research

The use of computations in mathematics educational research has been theorised in several ways. A widely used framework is computational thinking (CT), which positions programming as a disciplinary practice of computer scientists that provides transferable problem-solving skills across STEM fields (Wing, 2006; Weintrop et al., 2016). Research on CT in mathematics education has focused largely on how students leverage 'algorithmic' problem solving skills (see, e.g., Kotsopoulos et al. 2017; DeJarnette, 2019; Lockwood et al, 2020, 2021).

An alternative framing is computational literacy (CL) (diSessa 2000, 2018), which builds on Papert's (1980) view of computational media as expressive tools that allow learners to represent mathematical ideas as executable artefacts. In Papert's terms, computations "allow us to shift the boundary separating concrete and formal." (p. 21). CL takes a 'big picture' approach, focusing on how widespread literacy in expressive computational media reshapes how we learn, why we learn, and what counts as worthwhile knowledge.

Despite these, and other, theoretical developments, little is known about how mathematicians conceptualise the legitimacy of programming within their professional practices. Research in mathematics education focuses primarily on students' engagement with programming and the pedagogical affordances of computational tools, but has barely examined how disciplinary norms, identities, and epistemic commitments shape mathematicians' views of programming across research and teaching. This gap is striking given the longstanding ambivalence within the mathematical community toward computation and the growing institutional pressure to integrate programming into undergraduate curricula.

Taken together, the research reviewed in this section shows that mathematicians' views of programming in research and teaching are shaped by disciplinary norms and identities. However, little work has been done toward examining this link directly. Guided by a socio-cultural framework that we present in the next section, we adopt a "big-picture" view aligned with computational literacy: we consider programming not only as a technical skill, but as a practice that may be reshaping what mathematicians know, value, and do.

# Theoretical framework

## Cultural-Historical Activity Theory

We use Engeström's third-generation CHAT to conceptualise mathematicians' research and teaching as collective, object-oriented activities (Engeström, 2015). An activity is organised toward a shared 'object' – its driving purpose. For mathematicians, the object may be producing new mathematical knowledge, maintaining disciplinary standards, or educating students. These activities are historically enduring, partially overlapping, and can be a source of tension as people navigate across them. CHAT has been used in mathematics education research to study activity of teachers, students and their interactions (see, e.g., Núñez, 2009; Jaworski, 2010; Fredriksen & Hadjerrouit, 2020).

Engeström extends Vygotsky's classic mediational triangle by explicitly adding the collective factors that shape activity (see Figure 1). The 'Vygotskian' top part of the triangle represents how subjects act toward an object mediated by cultural tools. The nodes on the lower line (rules, community, and division of labour) capture what Engeström refers to as the hidden curriculum of an activity system (p.112, Engeström, 2015). This is the collective dimension of the activity: the norms that regulate action, the communities that confer legitimacy, and the distribution of roles and responsibilities that shape participation in the activity. These remind us that individual actions are always socially mediated.

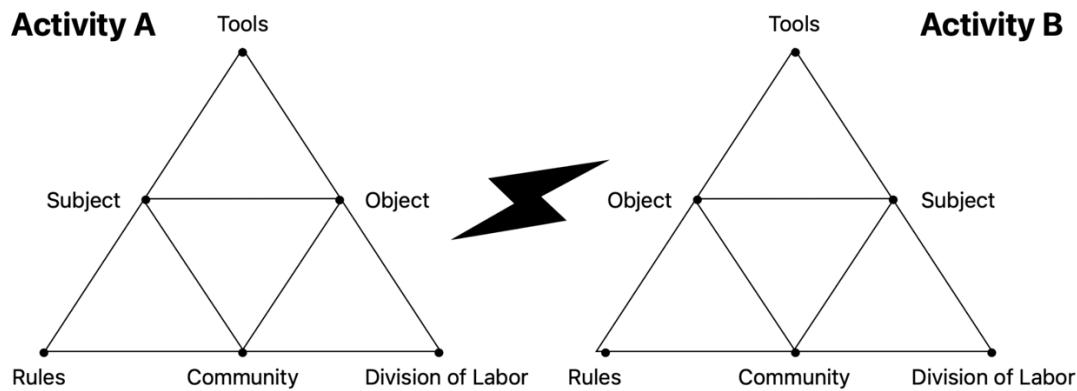

Figure 1: The unit of analysis of Engeström's third-generation CHAT is a pair of intertwined and interacting activity systems that have partially overlapping nodes.

Central to activity theory is the concept of *dialectic contradictions*. These are intrinsic contradictions, theorised as the driving forces of the activity. They are root causes of tension, conflict, change and innovation. As an example, Roth and Radford (2011) describe the fundamental contradiction of education: the object of learning is known to the teacher, but never to the student. This inevitably introduces an incommensurability of perspectives between teacher and student. Such contradictions are irremovable features of an activity, not problems to be solved.

To analyse tensions arising as mathematicians navigate between research and teaching, we draw on Engeström's (2015) notion of interacting activity systems. In Figure 1, we illustrate Engeström's unit of analysis, which is a pair of two interacting and partially overlapping activity systems. In our case, these correspond to the mathematicians' research and teaching activities.

Following Mwanza and Engeström (2005), we consider the nodes of the linked activity triangles as a checklist for questions to ask of our data set in order to identify underlying dialectic contradictions that drive the activities of teaching and research.

## Communities of practice

Communities of Practice (CoP) complements the systemic perspective of CHAT helping us understand how systemic tensions are taken up, reinterpreted, and enacted by individuals embedded in professional communities (Lave & Wenger, 1991; Wenger, 1998). We draw inspiration from Burton (2004) who used CoP to analyse 'discourse of power' in interviews with 70 mathematicians, and Jaworski (2010) who combined CoP with CHAT to analyse her own practices teaching mathematics at university level.

A central construct for our analysis from CoP is legitimate participation. Whether programming is positioned as full, peripheral, or marginal engagement in their practices reflects how different communities (and subcommunities) of mathematicians assign legitimacy to computational work. Burton (2004) shows how such judgments are rooted in epistemic and aesthetic values that shape what mathematicians recognise as authentic mathematical activity.

CoP foregrounds an engagement–reification duality, through which communities negotiate meaning by turning abstract ideas into shared objects. Terms such as "programming," "pure mathematics," or "mathematical work" are *reifications* whose meanings reflect and shape communal norms. Making sense of claims like "programming is not part of mathematics" requires attending to how such reified categories are produced, maintained, and contested.

Finally, CoP draws attention to boundaries as necessary consequences of specialisation. As mathematicians specialise, the meanings of reifications – and the norms of legitimate participation attached to them – changes. Moving across boundaries, such as between different subfields of mathematics, or between research and teaching, requires renegotiating these meanings, often surfacing tensions but also enabling new forms of practice. In our context, mathematicians navigate these boundaries and barriers as brokers, with "programming" functioning as a boundary object whose meaning shifts across settings.

## Research questions

In the following section, we operationalise CHAT's systemic search for contradictions through the human-centred lens of CoP to analyse how mathematicians describe the role of programming in their research and teaching.

Note that by 'computational tools' (or simply 'computations') we refer to ways of doing mathematical work that rely on literacy with expressive computational media (diSessa, 2018). Such work may or may not involve an actual computer; it can also take the form of reasoning shaped by familiarity with such media. By 'theoretical tools', we refer to ways of working that do not draw on such media. In line with our socio-cultural framing, we do not treat this distinction as objective, but as reflecting how computational media mediate what kinds of action participants view as possible.

This theoretical framing leads to the following research questions.

1. How does framing mathematicians' accounts in terms of their object orientation (theoretical vs. real-world) and tool orientation (theoretical vs. computational) help us make sense of coherent epistemic archetypes within an otherwise heterogeneous community?
2. Within each archetype, what epistemic rules govern the legitimate use of programming in research activity and teaching activity, respectively?
3. How do research and teaching differ in shaping mathematicians' views on the legitimacy of programming?

# Methodology

The aim of our study is to make sense of mathematicians' reflections on their use of programming in the contexts of their research and teaching. To do so, we conducted semi-structured interviews and analysed their accounts by drawing on CHAT and CoP (see, e.g., Jaworski, 2010). We take an interpretive stance, treating interviews as situated and co-constructed accounts of reflections of participants' epistemic orientations, identities, and interpretations of professional norms.

## Context and Participants

The study was conducted at the Faculty of Mathematics and Natural Sciences at the University of Oslo, which has a long-standing, institution-wide policy for integration of programming and scientific computation across STEM subjects for over 20 years (see Malthe-Sørenssen, 2016). This context allowed us to investigate mathematicians' views in an environment where programming is institutionally sanctioned and has received public acclaim.

Participants were recruited through purposeful, snowball sampling to capture diverse and contrasting perspectives on programming, and to ensure broad institutional representation. We interviewed fifteen mathematicians representing a wide range of career stages and specialisations (pure mathematics, applied mathematics, computational mathematics, and numerical analysis).

## Data Collection

Interviews were conducted individually between September and October 2023 and lasted 90 minutes to two hours. A semi-structured protocol guided each interview, organised around three areas:

- Research practices
- Teaching practices
- Institutional context

As interviews progressed, the protocol remained flexible, allowing participants to expand on themes emerging from previous conversations. All interviews were audio recorded, transcribed, and anonymised. Excerpts included in the paper were translated to English.

The interviews were conducted by the first author, an active research mathematician, which allowed participants to use professional jargon freely. This enabled added epistemic and technical depth when discussing proofs, modelling, or computational heuristics. This insider position, however, also carries risks of bias and over-alignment with participants' assumptions. To mitigate these risks, we approached the interviews as co-constructed dialogues in which the interviewer and participants jointly articulated meanings and interpretations. The second author served as an external critical voice during analysis and helped counterbalance potential insider bias.

## Analysis

Following Mwanza and Engeström (2005), we operationalise the nodes of Engeström's extended activity triangle as an analytic checklist that we systematically use to search for manifestations of contradictions in transcribed interview data. To ground this checklist in accounts of lived human experience, we draw on discourse indicators informed by CoP (Table 1).

| Table 1: The professional activity of 'doing mathematics' with illustrative quotes chosen from the data set. | |
|---|---|
| Theoretical construct | Focus of theoretical construct as analytic lens and operationalisation (using CoP) |
| Object | Collective purpose and goal of working as a mathematician. *"Mathematics is about exploring the thin membrane between the mathematically trivial and impossible."* |
| Subject | Identity and personal motive of why mathematicians do what they do. *"I just love theorems."* |
| Rules | Disciplinary norms and legitimacy of ways of how to do mathematics. *"I don't consider implementing algorithms to be part of my mathematical work."* |
| Tools | Material and immaterial tools used by mathematicians in their work. *"I use computers to check conjectures numerically. It allows me to do the theoretical work with more confidence and energy."* |
| Community | A map of the professional landscape of participation, legitimacy and boundaries. |
| | Table 1: The professional activity of 'doing mathematics' with illustrative quotes chosen from the data set. |
| | *"I find working on abstract problems without connection to reality to be boring and silly."* |
| Division of Labor | Use of agency to negotiate and navigate the mathematical community. *"I use computations to check conjectures, but most of my colleagues do not – even if they know that they'd probably get some new insight from doing this."* |
| Outcome | Tangible and symbolic products of professional mathematical practices. *"Almost all journals want […] an illustration based on real data."* |

Analysis proceeded iteratively. In an initial pass, the first author did a framework analysis using our operationalisation of CHAT to organise the data, followed by ongoing discussions within the author team in which interpretations, themes, and emerging hypotheses were continually challenged and refined. These interpretations were further stabilised through presentations at seminars and conferences (Olsen et al., 2025a, 2025b). Through this process, we developed an epistemic matrix capturing participants' orientations toward objects and tools (Table 2, below).

Our epistemic matrix guided a second phase of deductive coding, enabling focused comparison of object and tool orientations across research and teaching. Throughout an iterative process of analysis, the author team critically interrogated the evolving interpretations and theoretical commitments. The first author's position as a mathematician enabled triangulation with participants' publication records and research outputs when needed to clarify interpretations. To preserve anonymity and emphasise broader patterns, findings are reported at the level of archetypes rather than individuals.

# Findings

To address our research questions, we first present an epistemic matrix that emerged as part of our analysis. This matrix captures the dominant orientations mathematicians expressed toward the objects of their work and the tools they considered legitimate. We then use the matrix to describe how participants framed programming in their research and teaching, organised by the four archetypal practices represented in the matrix. Finally, we compare these accounts to examine differences between use of computational tools in research and teaching.

## The Epistemic Matrix: Four archetypal mathematical practices

Across interviews, participants differed systematically in their object (purpose) of research and the tools they regarded as legitimate for doing so. Through iterative analysis, we found that these differences could be organised along two intersecting dimensions:

- whether their object of research is oriented toward theoretical problems or problems that mix theoretical and real-world problems, and
- whether their preferred tool use is oriented toward theoretical tools or a mix between theoretical and computational tools.

As noted earlier, we follow diSessa (2018) in distinguishing theoretical tools from computational tools, where the latter are understood as literacies with expressive computational media. In referring to epistemic norms, we align with Weber et al. (2014), and mean the rules (in the CHAT sense) that govern what mathematicians seek to know and the means by which they justify that knowledge.

The resulting epistemic matrix (Table 2) helps us identify four archetypal ways of engaging in mathematics. Following our socio-cultural framework, we treat these archetypes as patterns of participation – recurrent ways in which mathematicians orient to objects and tools that, over time, become recognisable forms of practice.

Several mathematicians described engaging practices associated with multiple quadrants at different points in their career. However, we interpreted all but one to clearly articulate one orientation that was consistently more salient in how they described their daily work, their preferred problems, and their sense of what "counts" as mathematics. We therefore assigned the remaining mathematicians to the archetype best corresponding to the dominant points brought up in the interview. For narrative clarity we therefore refer to individual mathematicians as, say, "classical pure" or "computational applied", treating these labels as analytic constructs rather than self-ascribed identities. They reflect how participants positioned their practice through discourse, and not disciplinary categories.

| Table 2: Epistemic matrix (generic version) | | |
|---|---|---|
| | Theoretic object orientation | Pragmatic object orientation |
| Theoretic tool orientation | *"Classical pure mathematics"*<br><br>Prefers problems of theoretical interest and wants to pursue them with the aid of theoretical tools. If tools do not work, they tend to change problems.<br><br>*Frequency: 5 interviews* | *"Classical applied mathematics"*<br><br>Prefers problems of theoretical *and* real-world interest and wants to pursue them with the aid of theoretical tools. If tools do not work, they tend to change problems.<br><br>*Frequency: 2 interviews* |
| Pragmatic tool orientation | *"Computational pure mathematics"*<br><br>Prefers problems of theoretical interest, and wants to pursue them with theoretical and computational tools. If tools do not work, they are prepared to enhance the toolbox.<br><br>*Frequency: 3 interviews* | *"Computational applied mathematics"*<br><br>Prefers problems of theoretical *and* real-world interest, and wants to pursue them with theoretical and computational tools. If tools do not work, they are prepared to enhance the toolbox.<br><br>*Frequency: 5 interviews* |

## Mathematicians' views on the use of programming

We now use the epistemic matrix to guide our analysis of mathematicians' epistemic views on the use of programming in their research and teaching practices. These findings are summarised in Table 3.

|  | Theoretic object orientation | Pragmatic object orientation |
|---|---|---|
| Theoretic tool orientation | *"Classical pure mathematics"*<br>**Programming in research:**<br>Low or no engagement with computational tools. Do not see programming as mathematical work. Light-touch use is legitimate for 'pre-epistemic' work. Computational output is mostly absent from publications.<br>**Programming in teaching:**<br>Low engagement with programming in teaching. Substantive, theoretically framed integration is seen as legitimate, but should not disrupt analytic course identities. Light-touch support is accepted as helpful but has low epistemic visibility. Programming has low visibility in assessment. | *"Classical applied mathematics"*<br>**Programming in research:**<br>Skilled programmers who use computations as proof-of-concept for theoretical results. Do not see programming as mathematical work. Light-touch use is legitimate for 'pre-epistemic' work. Computational output has high visibility in publications.<br>**Programming in teaching:**<br>Low engagement with programming outside specialised methods courses. Legitimate integration requires substantive use and theoretical framing. Light-touch support is seen as helpful but has low epistemic visibility. Programming has low visibility in assessment. |
| Pragmatic tool orientation | *"Computational pure mathematics"*<br>**Programming in research:**<br>High but varied engagement with programming as legitimate mathematical work. All recount stories of being fascinated by the epistemic impact of computations. Computational artefacts have full legitimacy alongside theoretical results.<br>**Programming in teaching:**<br>Low engagement with programming in teaching, but with some notable exceptions. Legitimacy tied to perceived impact on understanding theoretical concepts. Light-touch support has low epistemic visibility. Programming has low visibility in assessment. | *"Computational applied mathematics"*<br>**Programming in research:**<br>High engagement with programming as a legitimate partner to analytic work on real-world problems. Computational artefacts have full legitimacy alongside theoretical results.<br>**Programming in teaching:**<br>Broad engagement with programming, but with some notable exceptions. Programming is seen to have intrinsic legitimacy, but tempered by perceived impact on understanding theoretical concepts. Programming has low visibility in assessment. |

Table 3: Epistemic matrix (research and teaching practices)

# Classical pure mathematicians

*Research practices*

Classical pure mathematicians describe their research as centred on abstract mathematical objects, which simultaneously serve as objects of study, tools of inquiry, and outcomes. As one senior mathematician describes, this tight alignment leaves little room for programming:

> "... there is very little in my research that uses programming or that uses any computer tools at all, really. [...] 'I tend to say that what is mathematically interesting is a very thin membrane between the trivial and the impossible — it is only this that you can answer with mathematics.' [...] You have to look for problems that you think are solvable with the tools you do have, after all. And that you can reach for in order to expand this toolbox a bit, right?"

While they report low engagement with programming, the classical pure mathematicians acknowledge the epistemic potential of the use of programming. As one said, "everyone agrees that having someone test such things [numerically] is a good thing. But that's not what is essential; it's the theorems you can prove..."

Although all had some familiarity with programming, and almost all had tried some exploratory numerical work, only one reported used it regularly as a 'light-touch' support for analytic work – occasionally checking conjectures and examples. He hesitated calling this programming, as he only ever used short snippets of code.

The output of programming is only visible in publications as support for conjectured results. Interestingly, they do not seem to talk much about computing with colleagues. Indeed, some were not able to comment on whether long-time collaborators program or not: "I'm sure he can program — he can do anything – but I don't think he bothers to do it."

*Teaching practices*

Classical pure mathematicians show low engagement with programming in their teaching, instead focusing on abstract mathematical objects studied using analytic tools such as theorems and proofs. They are generally positive to students learning how to program as part of their studies, but do not view courses in theoretical mathematics as the appropriate place to teach programming.

They talk about legitimate integration of programming into foundational courses as requiring substantive use and theoretical framing. One mathematician, who did not integrate computations in his courses, acknowledged that it is "amazing" how much can be learned from a simple for-loop. However, he was dismissive of a light-touch approach as the basis for meaningful integration: "I have been a bit embarrassed, because I feel that the programming we do here is [...] basically a loop with some setup beforehand and some cleanup afterwards. [...] It doesn't really build up any substantial programming competence." Instead, he talked about importing theory from numerical mathematics into classical theoretical courses as necessary to avoid 'trivial' use of computations.

However, this 'heavy-handed' approach was criticised by other participants as perturbing traditional course identities. Expressing the general view, one stated critically, "I believe slow learners, like me, [...] benefit from subjects that are like a comprehensible whole [...] instead of this crossover that only seems to help those who already understand everything."

Although they do not see it as substantive, the classical pure mathematicians speak positively about the use of a 'light-touch' of programming to support analytic learning. One notes that "weaker students" seemed to perform better after being "primed" by a little numerical exploration. Another described how using just the notion of programming had helped students connect ideas across disciplinary boundaries: "it was absolutely fantastic". Also those who did not engage with programming expressed that they were "basically positive" about having students "just sit and code a bit on the side, you know, simply to develop intuition." At least if they already knew how to program. However, as one added, "it maybe should be enough to [let students] just sit there with pen and paper".

Across the archetype, programming is nearly invisible in assessment. One says, "but we expect them to read some simple Matlab output. But you can solve the problems without being able to, but then you may have to do a bit of computing [by hand]."

## Classical applied mathematicians

*Research practices*

Classical applied mathematicians describe their research practice as predominantly oriented toward theoretical work on numerical algorithms for mathematical models having some real-world relevance. Although they are skilled programmers who use computations overtly in their work, they hesitate to describe it as part of their mathematical practices. As one explains:

> "The research work consists of developing numerical methods and then analyzing those methods using functional analysis and other tools. And the programming of these methods is not the most important [...] It is mostly to illustrate that this is not loose talk. [...] I consider the theory for and analysis of the numerical methods as the truly mathematical part."

While some of the classical applied mathematicians first played down the epistemic role of programming, they expressed more nuance when discussing the use of programming in more detail. They all report using programming as light-touch epistemic support for analytic work, though, as one notes, this happens "not so often." At the same time, they emphasise that computations serve a deeper purpose than merely cosmetic illustration – indeed, omitting them can appear suspicious: "it would seem to me a sign that there is something they are trying to hide, that this doesn't actually work in practice."

Interestingly, the explicitly stated norm of publications having to be centered on theory accompanied by computational proof-of-concept had just as clearly stated exceptions. At the one extreme, results about numerical methods that could not easily (or at all) be implemented were excused from needing computational support. As one said, "[sometimes] we can't be bothered to look at the fully discrete case". At the other extreme, they consider it legitimate to publish papers based solely on computational results when working on physical models where little or no theory was known – what they describe as a "physicist's approach." In this way, they voice how their practice shares a boundary with that of the classical pure and computational applied mathematicians.

*Teaching practices*

Classical applied mathematicians show the same level of low engagement with programming in foundational courses as the classical pure mathematicians. They find it difficult to teach programming, and prefer to compartmentalise this to specialised methods courses.

They see legitimate integration of programming in their study programs as requiring theoretical framing in numerical theory, yet are critical of efforts to integrate programming by moving theory from numerical mathematics into other courses as arbitrary: "It seems to me that someone who has no clue about programming or numerical analysis has been asked to insert some programming into a course [...] It feels a bit tacked on."

When asked about more light handed approaches, some express skepticism: "I believe that you really must understand mathematical theory first before you can do programming. [...] It certainly isn't easy to learn it the other way around… to somehow rely more on mathematical computation on a computer and understand the theory." This echoes their view from research, that programming is primarily a tool to be used to illustrate already obtained mathematical results, and not as a tool for learning mathematics.

As with the classical pure mathematicians, the applied mathematicians marginalise the role of programming in assessment, even talking about having let students off the hook if they were able to do the analytic part of the course work.

# Computational pure mathematicians

*Research practices*

Computational pure mathematicians describe their practice as centered on abstract mathematical objects, which they pursue by drawing on a combination of analytic and computational tools. They are the most heterogeneous group, but linked by stories they tell of becoming fascinated by how computational tools allowed them to gain insight on pure objects far beyond what was possible with pen and paper computations. As one says, "being able to use these specialized programs to create, compute, and understand more examples makes a huge difference for me."

Computational pure mathematicians leverage programming in several ways, ranging from light-touch numerical exploration to deep algorithmic implementation. Some use computations much like their classical pure colleagues – checking conjectures or generating examples – but do so openly, providing symbolic or numeric scripts alongside publications. As one says, matter of factly, "you should always explain what you do." Others take programming further, implementing algorithms to probe the structure of abstract objects and treating computation as a complementary and legitimate source of insight.

They go beyond the light computational support of the classical pure mathematicians, in that they are highly skilled programmers, able to leverage supercomputers if need be. A few describe computation as offering fundamentally new ways of knowing, arguing that developments such as generative AI reveal limits of purely analytic methods and make combined analytic-computational approaches indispensable. If not, one says, theoretical mathematics risks being reduced to "philosophy".

*Teaching practices*

While the computational pure mathematicians hold a high view of the legitimacy of the use of computations to learn mathematics, they show uneven engagement with programming in teaching. Interestingly, they point to pedagogical concerns for the learning of analytic mathematics as reasons both for high and low engagement with programming in their teaching.

They see the legitimacy of the use of programming in their teaching to be tied to what they see as necessary for learning about analytic concepts. For instance, one of them justifies low engagement with programming in a foundation course by noting that "you can do the questions by hand, so I didn't tell them [to use programming]". In more advanced courses, he describes how we demonstrate computational tools that students can use to computationally explore examples that are hard, or intractable, by purely analytic means.

This also indicates that, as with the other archetypes, light-touch use of programming to support learning of analytic concepts has low visibility. Further evidence of this is provided when another describes observing how students who struggle analytically, struggle even more when also asked to use computational tools. For this reason, he sees his role as offering courses where students can focus on analytic material. He says, "my teaching has always [...] focused on the theoretical part. It's not to [exclude] the other side, but to be complementary to the computations. I think it's important that they get both parts." However, some argue that a maximalist stance to the integration of computations could be beneficial to students struggling with theory: "I think… for someone who is quite strong at programming, regarding their ability to learn the theory, I think it helps. They see it more easily; they understand [the theory] better through code than through mathematical proof." However, he stops short of implementing this philosophy. He describes rarely bringing up programming in lectures, leaving it mostly as a complementary perspective in course material. He says, "You should not expose them to too much code in all of [our] mathematics courses, but it is beneficial for them."

Across the archetype, programming is nearly invisible in assessment. One says, "the only place we force them to do any programming is on homework assignments."

# Computational applied mathematicians

*Research practices*

Computational applied mathematicians describe their research practice as developing and implementing mathematical models of real-world phenomena through an interplay of analytic tools, computational implementation and collaboration with external partners. They come across as skilled programmers who are genuinely fascinated by working in the interface between real-world data, computations and theoretical mathematics.

> "It's exactly at that intersection that I find it enjoyable to work, so to speak — where you can have a kind of geometric intuition for how the problem looks [...] and then, from that understanding, use the right numerical moves to carry it far enough that you can find a concrete answer in a given situation."

They position computations and analytic tools as forming a legitimate partnership where each supports the other. As one says, "for me there isn't a clear separation in the sense that some problems require analysis, others modelling, others require programming." Further illustrating this point, one says about colleagues who did theoretical work on numerical algorithms, but without having computations in mind: "You could say that he produced formulas that were, insofar as they go, completely unassailable, but which were almost impossible to compute with."

Programming has high visibility in the outcome of their work, which even includes producing and maintaining software on a commercial level.

*Teaching practices*

Overall, computational applied mathematicians show a much higher engagement with programming in their teaching than other archetypes, most integrating computations to some degree. However, their engagement with computations seems to be tempered by concern for students' analytic abilities. One of them explains why he does not use any programming in a foundational course he teaches: "the challenge for people to get through the first-year course I teach does not lie in their ability to use a computer. Rather, it lies in their ability to understand [analytic mathematics]."

Almost all are quick to describe various ways in which they use light-touches of programming in support of learning theory: "there are some for-loops; it is not an advanced type of programming. We use it both so that they actually produce results, and because a bit of programming helps them understand the methods better." However, some are frustrated as they see such light-touch support as all that students get experience doing, and argues that there is a need for more substantial computational projects in their study program: "I mean, anyone can do that. If you are going to do computations properly you have to do something beyond [using snippets of code]. Something like… you need to have gotten your hands dirty programming a bit …"

Across the archetype, programming is nearly invisible in assessment. One says "we think of programming and coding more as an aid. It is not something we see as essential to test."

## Comparison of research and teaching practices

Across all four archetypes, mathematicians seem to carry the same underlying epistemic rules from research into teaching, exhibiting epistemic continuity across settings. Pure mathematicians maintain 'analytic exclusivity'; applied mathematicians uphold analytic exclusivity with pragmatic concessions; computational applied mathematicians rely on analytic–computational partnership; and computational pure mathematicians treat computation as legitimately analytic. These rules anchor their research practices and remain intact when they shift to teaching.

At the same time, the object-outcome shift from research to teaching, described by Burton (2004), reshapes how programming is used. Teaching focuses on the established outcomes of mathematics – definitions, theorems, and analytic techniques – rather than on the practices by which such knowledge is produced. Because foundational courses foreground analytic content and conventional assessment formats, the epistemic role of computation tends to lose visibility. Even mathematicians who rely heavily on programming in research restrict its use in teaching, citing students' weaker analytic preparation and the need to preserve course identity.

These dynamics give rise to a recurring tension between heavy-handed and light-touch approaches to integrating programming. Heavy-handed approaches seek theoretically legitimate integration – often by importing numerical analysis material or programming instruction into pure courses – but this is seen to disrupt course identities and disadvantage students without strong analytic ability. By contrast, light-touch computational use (e.g., exploratory scripts, numerical checks, visualisations, or algorithmic framings without actual code) is widely viewed as pedagogically helpful, especially for students who struggle analytically. Yet such practices remain marginalised because they are rarely recognised as "real" programming and therefore have low epistemic visibility.

Finally, programming remains largely invisible in assessment. Across archetypes, mathematicians avoid examining programming, describing it as difficult to grade or peripheral to "essential" mathematical competence. Even those who rely heavily on computation in research sometimes overlook incomplete computational work if students demonstrate analytic proficiency. We interpret this as a strong signal of low epistemic legitimacy of computations in teaching, across all archetypes, and it contrasts sharply with the high visibility of computational contributions in many participants' research publications.

# Discussion

## Synopsis

In this study, we leveraged a CHAT and Communities of Practice to analyse interviews with 15 mathematicians with respect to how they articulate the legitimacy of programming in their research and teaching. We now discuss the analysis in light of the three research questions we formulated at the end of the Theoretical framework section.

RQ1: *How does framing mathematicians' accounts in terms of their object orientation (theoretical vs. real-world) and tool orientation (theoretical vs. computational) help us make sense of coherent epistemic archetypes within an otherwise heterogeneous community?*

Despite longstanding institutional support, substantial variation persisted in how mathematicians positioned programming as an epistemic tool. Our analysis showed that this variation in participants' discourse was not purely idiosyncratic, but patterned. To see this, we organised mathematicians according to object and tool orientation, identifying four epistemic archetypes: classical pure, classical applied, computational applied, and computational pure. Within each archetype, mathematicians expressed internally coherent views about what counted as meaningful problems, acceptable tools, and legitimate outcomes, enabling us to reduce the substantial heterogeneity described in prior literature. We saw that the pure applied and computational pure mathematicians expressed patterns of boundary crossing more than the classical pure and computational applied mathematicians.

*RQ2: Within each archetype, what epistemic rules govern the legitimate use of programming in research activity and teaching activity, respectively?*

In research, all archetypes abided by some form of coordination between the use of analytic and computational tools – from strict analytical exclusivity to analytic–computational partnership. These orientations shaped whether programming was framed as exploratory 'pre-epistemic' scaffolding, as numerical confirmation, or as a tool for establishing legitimate insight. In teaching, mathematicians' epistemic rules largely carried over from research; however, because teaching in foundational courses foregrounds outcomes of established knowledge rather than the processes of inquiry that shape research, computational work had markedly lower epistemic visibility. This object–outcome shift systematically constrained the role programming could play in teaching – even among those who depend on computation professionally.

*RQ3: How do research and teaching differ in shaping mathematicians' views on the legitimacy of programming?*

A central finding was the tension between heavy-handed and light-touch approaches to computation. Heavy-handed integration, based on adding theory from numerical analysis into traditionally purely analytic courses, offered epistemic legitimacy but disrupted established course identities and often privileged analytically strong students. Light-touch programming – small snippets used to support conceptual understanding – was more pedagogically accessible but was not recognised by most mathematicians as "real programming," leaving it epistemically marginalised. Across archetypes, programming remained largely absent from assessment: even mathematicians who rely heavily on computation in research avoided including it in assessment, indicating a conflicted view of epistemic status in teaching.

# Contributions and implications for education

We offer a new perspective on the status gap suggested by Broley et al. (2018) between mathematicians' use of programming in research and teaching: it reflects not only cultural attitudes or institutional constraints, but deeper epistemic continuities rooted in mathematicians' preferred modes of reasoning. Mathematicians carry their research-based norms into teaching, but because teaching foregrounds the outcomes of mathematics and relies on analytic assessment, programming loses epistemic visibility.

Similarly, our results suggest a reframing of Burton's (2004) explanation of the research–teaching gap. Like Burton, we find that the object of teaching is taken to be the outcomes of research, marginalising exploratory practices central to research. But while Burton attributes the research-teaching gap in practices to mathematicians enacting different epistemologies (in her more broad, sociological, sense) our findings suggest that when considering epistemologies in the more narrow sense of Weber et al. (2014), the gap can be understood by mathematicians displaying epistemic continuity across practices in combination with the object-outcome shift. In light of the methodological challenges identified by Weber et al. (2020), we suggest that analytic framings which foreground such continuities may offer a productive path for future research on mathematicians' practices.

Furthermore, our use of an epistemic matrix helped us handle some of the methodological issues raised by Weber et al. (2020) by characterising mathematicians through their object and tool orientations, and observing how concepts such as 'programming' are reified differently across boundaries of practice. This suggests that epistemic orientation may be a productive analytic lens for future research on mathematicians' professional and teaching practices.

In light of how mathematicians articulate legitimate uses of programming in their professional work, a practical implication for education is to explicitly legitimize light-touch epistemic programming. Making its purposes visible, articulating how it supports analytic understanding, and aligning modest assessment with these aims could help narrow the gap between mathematicians' research practices and their teaching.

# Limitations and future research

Consistent with our socio-cultural and interpretivist stance, we see our results as context-dependent co-constructions between the team of authors and participants.

The insider position of the interviewer, a research mathematician, enabled deep technical dialogue but also introduced a risk of over-alignment with participants' assumptions. We mitigated this by an extended period of analysis where partial results were shared and discussed with mathematics and physics education researchers, as well as professional mathematicians, promoting interpretive triangulation.

We acknowledge that our data was collected in fall 2023, and analysed over the following two year period, and are inevitably colored by a moment in time when large language models were rapidly evolving and starting to impact programming practices. Rather than a limitation, this suggests an opportunity for future research: several mathematicians were already positioning AI tools epistemically, offering a natural extension of our framework.

Our findings point to several directions for further investigation. First, it would be of interest to use the same analytic approach at other sites to shed light on how our findings extend beyond the specific context of this study and add further nuance. Second, our results raise questions about how disciplinary boundaries shape the meanings attached to programming, and how boundary-crossers – such as classical applied and computational pure mathematicians – may facilitate innovation and change. Third, in this paper, we foregrounded an analysis on epistemic norms and beliefs. In future work, it would be interesting to explore, in more detail, how and why programming is used by the interviewed mathematicians in their research and teaching. Fourth, while institutional factors appear to mediate mathematicians' decisions to adopt programming, mathematicians appear to have a large degree of operational agency. A more detailed understanding of the interaction between institutional pressure and individual agency would be interesting for efforts to enact reforms involving the use of computations. Fifth, the epistemic matrix helped us interpret mathematicians' orientations toward programming as a tool for learning. Future research could investigate whether it also offers insight into how mathematicians position generative AI epistemically (cf., Fraser et al., 2024a, 2024b).

**Conflict of interest:** On behalf of all authors, the corresponding author states that there is no conflict of interest.

# References

Broley, L., Caron, F., & Saint-Aubin, Y. (2018). Levels of programming in mathematical research and university mathematics education. *International Journal of Research in Undergraduate Mathematics Education, 4,* 38-55. http://dx.doi.org/10.1007/s40753-017-0066-1

Burton, L. L. (2004). *Mathematicians as enquirers: Learning about learning mathematics (Vol. 34).* Springer Science & Business Media.

Buteau, C., Jarvis, D., & Lavicza, Z. (2014). On the integration of computer algebra systems (CAS) by Canadian mathematicians: Results of a national survey. Canadian Journal of Science, Mathematics, & Technology. Education, 14(1), 1–23. http://dx.doi.org/10.1080/14926156.2014.874614

Darragh, L. (2022). Brokering across the divide: Perspectives of mathematicians involved in education. *The Journal of Mathematical Behavior, 67*, 100989. https://doi.org/10.1016/j.jmathb.2022.100989

DeJarnette, A. F. (2019). Students' challenges with symbols and diagrams when using a programming environment in mathematics. *Digital Experiences in Mathematics Education, 5*(1), 36–58. https://doi.org/10.1007/s40751-018-0044-5

DiSessa, A. A. (2000). *Changing minds: Computers, learning, and literacy*. Mit Press.

diSessa, A. A. (2018). Computational literacy and "the big picture" concerning computers in mathematics education. *Mathematical thinking and learning, 20(1)*, 3–31. http://dx.doi.org/10.1080/10986065.2018.1403544

Engeström, Y. (2015). *Learning by Expanding: An Activity-Theoretical Approach to Developmental Research.* Cambridge University Press.

Fraser, M., Granville, A., Harris, M. H., McLarty, C., Riehl, E., & Venkatesh, A. (2024a). Will machines change mathematics? *Bulletin (New Series) of the American Mathematical Society, 61(2),* 201–202. https://doi.org/10.1090/bull/1836

Fraser, M., Granville, A., Harris, M. H., McLarty, C., Riehl, E., & Venkatesh, A. (2024b). Will machines change mathematics? *Bulletin (New Series) of the American Mathematical Society, 61(3),* 373–374. https://doi.org/10.1090/bull/1842

Helge Fredriksen & Said Hadjerrouit (2020) An activity theory perspective on contradictions in flipped mathematics classrooms at the university level*, International Journal of Mathematical Education in Science and Technology, 51(4)*, 520-541, https://doi.org/10.1080/0020739X.2019.1591533

Hamami, Y., & Morris, R. L. (2020). Philosophy of mathematical practice: a primer for mathematics educators. *ZDM, 52(6)*, 1113-1126. https://doi.org/10.1007/s11858-020-01159-5

Hersh, R. (1997). *What is mathematics, really?* Oxford University Press.

Jaworski, B. (2010). The practice of (university) mathematics teaching: Mediational inquiry in a community of practice or an activity system. *CERME 6–WORKING GROUP 9,* 1585–1594.

Kotsopoulos, D., Floyd, L., Khan, S. et al. (2017). A Pedagogical Framework for Computational Thinking. *Digit Exp Math Educ* **3**, 154–171 https://doi.org/10.1007/s40751-017-0031-2

Lave, J., & Wenger, E. (1991). *Situated learning: Legitimate peripheral participation*. Cambridge university press.

Lockwood, E., DeJarnette, A. F., & Thomas, M. (2019). Computing as a mathematical disciplinary practice. *Journal of Mathematical Behavior, 54*. https://doi.org/10.1016/j.jmathb.2019.01.004

Lockwood, E., De Chenne, A. Enriching Students' Combinatorial Reasoning through the Use of Loops and Conditional Statements in Python. *Int. J. Res. Undergrad. Math. Ed.* **6**, 303–346 (2020). https://doi.org/10.1007/s40753-019-00108-2

Lockwood, E., & Mørken, K. (2021). A call for research that explores relationships between computing and mathematical thinking and activity in RUME. *International Journal of Research in Undergraduate Mathematics Education, 7(3),* 404-416. https://doi.org/10.1007/s40753-020-00129-2

Lovelace, A. (n.d.). Computer, enhance. *Lapham's Quarterly.* Retrieved April 5, 2025, from https://www.laphamsquarterly.org/technology/computer-enhance

Malthe-Sørenssen, A., Hjorth-Jensen, M., Langtangen, H. P., & Mørken, K. (2016). Integrating computation in the teaching of physics. *Uniped, 38*. Retrieved from http://hplgit.github.io/cse-physics/doc/pub/uniped15.html

Moschkovich, J. (2007). Examining mathematical discourse practices. *For the Learning of Mathematics,* 27(1), 24–30.

Mwanza, D., & Engeström, Y. (2005). Managing content in e-learning environments. *British Journal of Educational Technology, 36(3),* 453–463. https://doi.org/10.1111/j.1467-8535.2005.00479.x

Nardi, E. (2008). *Amongst mathematicians: Teaching and learning mathematics at university level (Vol. 3).* Springer Science & Business Media.

Núñez, I. (2009). Contradictions as Sources of Change: A literature review on Activity Theory and the Utilisation of the Activity System in Mathematics Education. *Educate, 9(3),* 7-20.

Olsen, J. F., Lockwood, E., & Odden, T. O. B. (2025a). *Between the trivial and the impossible: CHAT as an analytic lens to study the role of programming in mathematical research.* In Proceedings of the Fourteenth Congress of the European Society for Research in Mathematics Education (CERME14) (No. 20).

Olsen, J. F., Odden, T. O. B., & Lockwood, E. (2025b). *Exploring Discrepancies in Mathematicians' Use of Programming in Teaching and Research.* In S. Cook, B.P. Katz, & K. Melhuish (Eds.), Proceedings of the 27th Annual Conference on Research in Undergraduate Mathematics Education (pp. 1102–1107). SIGMAA on RUME. Alexandria, VA.

Papert, S. (1980). *Mindstorms: Children, Computers, and Powerful Ideas*. Basic Books: New York, 1980. 230 pages.

Polya, G. (2014). *How to solve it*. Princeton University Press.

Roth, W. M., & Radford, L. (2011). *A cultural-historical perspective on mathematics teaching and learning (Vol. 2)*. Springer science & business media.

Rasmussen, C., Wawro, M., & Zandieh, M. (2015). Examining individual and collective mathematical progress. *Educational Studies in Mathematics, 88*, 259–281. http://dx.doi.org/10.1007/s10649-014-9583-x

Weber, K., Dawkins, P., & Mejía-Ramos, J. P. (2020). The relationship between mathematical practice and mathematics pedagogy in mathematics education research. *ZDM, 52(6)*, 1063-1074. http://dx.doi.org/10.1007/s11858-020-01173-7

Weber, K., Inglis, M., & Mejia-Ramos, J. P. (2014). How mathematicians obtain conviction: Implications for mathematics instruction and research on epistemic cognition. *Educational Psychologist*, *49*(1), 36-58. https://doi.org/10.1080/00461520.2013.865527

Weintrop, D., Beheshti, E., Horn, M., Orton, K., Jona, K., Trouille, L., & Wilensky, U. (2016). Defining computational thinking for mathematics and science classrooms. *Journal of science education and technology, 25(1)*, 127-147. http://dx.doi.org/10.1007/s10956-015-9581-5

Wenger, E. (1998). *Communities of Practice: Learning, Meaning and Identity*, Cambridge UK: Cambridge University Press.

Wing, J. M. (2006). Computational thinking. *Communications of the ACM, 49(3)*, 33-35. http://dx.doi.org/10.1145/1118178.1118215